\documentclass[12pt,twoside]{amsart}

\usepackage{amsfonts}
\usepackage{amsmath}
\usepackage{latexsym}

\newtheorem{thm}{{\sc Theorem}}[section]

\newtheorem{lem}[thm]{{\sc Lemma}}
\newtheorem{cor}[thm]{{\sc Corollary}}
\newtheorem{defn}[thm]{{\sc Definition}}
\newtheorem{claim}[thm]{{\sc Claim}}

\newcommand{\qqed}{\hspace*{\fill} $\Box$}

\newcommand{\quadric}{q_{Z}}
\newcommand{\quadrics}{q_{{\mathcal Z}_{s}}}
\newcommand{\quadricsurface}{Q}
\newcommand{\cupsurface}{S}
\newcommand{\PP}{{\mathbb P}}
\newcommand{\CC}{{\mathbb C}}
\newcommand{\twoform}{\omega}
\newcommand{\pulltwoform}{\tilde{\omega}}
\newcommand{\pulltwoforms}{\tilde{\omega}_{s}}
\newcommand{\barpulltwoforms}{\bar{\tilde{\omega}}_{s}}

\title[Fujiki Relation]{Fujiki relation on symplectic varieties}
\author{Daisuke Matsushita}
\date{}
\address{Division of Mathematics, Graduate School of Science,
         Hokkaido University,  Sapporo, 060-0810 Japan}
\email{matusita@math.sci.hokudai.ac.jp}

\begin{document}
\maketitle

\begin{abstract}
 We generalize Fujiki relation
 of Beauville-Bogomolov quadratic form on a projective symplectic
 variety. As an application, we study a fibre space structure of 
 a projective symplectic variety. 
\end{abstract}

\section{Introduction}
 We start with the definition of a {\itshape symplectic variety}.
\begin{defn}
 A compact K\"{a}hler variety $Z$ is said to be
 a symplectic variety if $Z$ satisfies the following 
  conditions:
       For a smooth locus $U$ on $Z$, there exists a
       nondegenerate holomorphic $2$-form $\twoform$ on $U$.
       This form is extended to a regular form $\pulltwoform$
       on $\tilde{Z}$,
       where  $\nu : \tilde{Z} \to Z$
       is a resolution of singularities of $Z$.
\end{defn}

\noindent
 In \cite[Theorem 8 (2)]{namikawa1}, Namikawa
 induced a quadratic form $\quadric$
 on $H^{2}(Z, \CC)$, which is a natural extension
 of Beauville-Bogomolov quadratic form 
 defined in \cite[Th\'{e}or\`{e}me 5]{beauville2} 
 on singular varieties.
 In \cite[Theorem 4.7]{fujiki}, Fujiki proved that
 Beauville-Bogomolov quadratic form has special relation
 with cup products.
 We prove that $\quadric$ has same properties.
\begin{thm}\label{main}
 Let $Z$ be a $2n$-dimensional projective symplectic
 variety and
 $\nu : \tilde{Z} \to Z$ a resolution of singularities
 of $Z$. Assume that
\begin{enumerate}
 \item The codimension of the singular locus of $Z$ is greater than
       four.
 \item $Z$ has only ${\mathbb Q}$-factorial singularities.
 \item $\dim H^{1}(Z,{\mathcal O}_{Z}) = 0$ and
       $\dim H^{2}(Z,{\mathcal O}_{Z}) = 1$.
\end{enumerate}
 According to \cite[Theorem 8 (2)]{namikawa1},
 we define
 the quadratic form $\quadric$ on $H^{2}(Z,{\mathbb C})$ by
$$ 
 \quadric (\alpha ):= n\int_{\tilde{Z}} 
 (\pulltwoform \bar{\pulltwoform})^{n-1}\tilde{\alpha}^2 +
 (1-2n)\left(\int_{\tilde{Z}} 
          \pulltwoform^{n}\bar{\pulltwoform}^{n-1}\tilde{\alpha}
            \int_{\tilde{Z}} 
          \pulltwoform^{n-1}\bar{\pulltwoform}^{n}\tilde{\alpha}
      \right),
$$
 where 
 $\tilde{\alpha} := \nu^{*}\alpha$ and $\alpha \in H^{2}(Z, \CC)$.
 Then $\quadric (\alpha)$ satisfies the following equation:
$$
  C_{Z} \quadric (\alpha )^{n} = \alpha^{2n}.
$$
 Note that a constant $C_{Z}$ depends on only $Z$.
 Moreover the index of $\quadric$ is $(3,h^{2}(Z, \CC ) - 3 )$.
\end{thm}
\noindent
{\scshape Remark 1. \quad} By \cite[Theorem]{namikawa2},
 when $Z$ is symplectic
  varieties with only terminal singularities, the condition
 of codimension of  singular locus is always satisfied.

\noindent
{\scshape Remark 2. \quad} Namikawa obtained
 the index of $\quadric $ in \cite[Corollary 8]{namikawa3}
 by different method.

\begin{thm}\label{main2}
 Let $Z$ be a projective symplectic variety which satisfies
 the conditions of Theorem \ref{main} and $D$
 a Cartier divisor on $Z$. Then 
 Riemann-Roch formula of $D$ is expressed as follows:
$$
 \chi (D) = \sum_{k=0}^{n} a_{k}\quadric (D)^{k},
$$
 where $a_{k}$ $(0 \le k \le n)$ are constants which depend on only $Z$.
\end{thm}
 
\noindent
 As an application of Theorem \ref{main}, we obtain the following
 result:
\begin{cor}\label{application}
 Let $Z$ be a $2n$-dimensional
 projective
 symplectic variety which satisfies the conditions of Theorem \ref{main}.
 Assume that there exists a surjective morphism $f : Z \to B$
 over a projective normal variety $B$ and a general fibre $F$ of $f$
 has positive dimension.
 Then $F$ and $B$ has the following properties:
\begin{enumerate}
 \item  $\dim B = n$ and Picard number of $B$ is one.
 \item  $F$ is an $n$-dimensional abelian variety.
 \item  For the singular locus $Z_{sing}$ of $Z$,
        $f(Z_{sing})$ forms a proper closed subset of $B$ and
        the restriction $\omega$ to $F$ is identically zero.
 \item  For every effective divisor $D$ of $Z$,
        $\dim f(D) \ge \dim B -1$.
\end{enumerate}
\end{cor}
\noindent
{\scshape Remark 3. \quad}
 Comparing the above corollary with \cite[Theorem 2]{matsu},
 there exist two difference: One is ampleness of $-K_B$ and
 another one is ${\mathbb Q}$-factoriality of $B$.
 In section 3, we construct
 an example such that $-K_B$ is not ample. 

\noindent
 We prove Theorem \ref{main} and \ref{main2} in Section 2.
 A proof of Corollary \ref{application} and
 an example are given in Section 3.

\noindent
{\scshape Acknowledgement. \quad} The author express his thanks to
 Professors A.~Fujiki and Yo.~Namikawa for their
 advice and encouragement.

\section{Proof of Theorem \ref{main}}
 \noindent
(2.1) We start with the investigation of the rank of $\quadric$.
\begin{lem}\label{rank}
 The rank of $\quadric$ is greater than three.
\end{lem}
\noindent
{\scshape Proof. \quad}
 Let $\nu : \tilde{Z} \to Z$ be a resolution of singularities of $Z$ and
 $\pulltwoform$ a extension of $\twoform$ on $\tilde{Z}$.
 By the definition of $\quadric$,
 $\quadric (\pulltwoform + \bar{\pulltwoform}) >0$ and
 $\quadric (\pulltwoform - i\bar{\pulltwoform}) > 0$.
 Let $\alpha$ be an ample divisor on $Z$.
 If we choose a suitable exceptional divisor $E$, $\nu^{*}\alpha - E$
 becomes an ample divisor. From \cite[Remark (1) (p24)]{namikawa1},
 $\pulltwoform^{n-1} |_{E} = 0$. Hence
$$
 \quadric (\alpha) = n\int_{\tilde{Z}} 
 (\pulltwoform \bar{\pulltwoform})^{n-1} (\nu^{*}\alpha - E)^{2}.
$$
 Thus $\quadric (\alpha) >0$ by Hodge-Riemann bilinear relation.
\qqed

\noindent
(2.2) Let ${\mathcal D}$ be the Kuranishi space of $Z$,
 ${\mathcal Z}$ the Kuranishi family and
 ${\mathcal Z}_{s}$ the fibre at $s \in {\mathcal D}$.
 By \cite[Theorem 8]{namikawa1}, 
 for every point of $s \in {\mathcal D}$,
 there exists an isomorphism 
 $\phi_{s} :  H^2(Z,\CC ) \to H^2({\mathcal Z}_{s},\CC ) $. 
 The period map $p : {\mathcal D} \to\PP (H^2(Z,\CC ))$ defined
 by $p(s) := \phi^{-1}_{s} (\twoform_{s})$,
 where $\twoform_{s}$ is a symplectic form of ${\mathcal Z}_{s}$.
\begin{lem}\label{prep1}
 Let $\nu_{s} : \tilde{{\mathcal Z}_{s}} \to
{\mathcal Z}_{s}$ be a resolution of 
${\mathcal Z}_{s}$.
 We define a quadric form on $H^{2}({\mathcal Z}_{s}, \CC )$ by
$$
 \quadrics (\alpha) := 
 n\int_{\tilde{{\mathcal Z}_{s}}} 
 (\pulltwoforms \barpulltwoforms)^{n-1}\tilde{\alpha}^2 +
 (1-2n)\left(\int_{\tilde{{\mathcal Z}_{s}}} 
          \pulltwoforms^{n}\barpulltwoforms^{n-1}\tilde{\alpha}
            \int_{\tilde{{\mathcal Z}_{s}}} 
          \pulltwoforms^{n-1}\barpulltwoforms^{n}\tilde{\alpha}
      \right),
$$
 where $\tilde{\alpha} := \nu^{*}_{s}\alpha$ and
 $\alpha \in H^{2}({\mathcal Z}_{s}, \CC )$.
 Then $\quadrics$ is not depend the choice of
 $\tilde{{\mathcal Z}_{s}}$ and 
$$
 \quadric (\alpha) =
 \quadric (\phi^{-1}_{s}(\twoform_{s} + \bar{\twoform}_{s}))
  \quadrics (\phi_{s}(\alpha)).  
$$
\end{lem}
\noindent
{\scshape Proof. \quad} 
 We prove that
 $\quadrics$ is independent of the choice of resolution
 by similar argument in the proof of \cite[Theorem 8 (2)]{namikawa1} and
 obtain the rank of $\quadrics$ is greater than
 three by similar argument of Lemma \ref{rank}.
 The period map 
 $p_{s} : {\mathcal D} \to \PP (H^{2}({\mathcal Z}_{s},\CC ))$
 is defined
 by $p_{s}(t) := \phi_{s}(p(t))$. 
\begin{claim}
 Let $ \quadricsurface_{s} := \{  \quadrics (\alpha ) = 0 |
     \alpha \in  H^{2}({\mathcal Z}_{s}, \CC) \}$.  
 Then $p_{s}({\mathcal D}) \subset \quadricsurface_{s} $.
\end{claim}
\noindent
{\scshape Proof. \quad} 
 Let $\alpha \in H^{2}({\mathcal Z}_{s}, \CC )$. 
 For the smooth locus of ${\mathcal U}_{s}$ of ${\mathcal Z}_{s}$,
 there exists Hodge decomposition of $H^{2}({\mathcal U}_{s},\CC )$
 and an isomorphism 
 $H^{2}({\mathcal U}_{s}, \CC ) \cong H^{2}({\mathcal Z}, \CC )$ by
 \cite[Theorem 8]{namikawa1}. Hence
 we write $\alpha$ as 
 $\alpha = a\twoform_{s} + w + b\bar{\twoform}_{s}$.
 By direct calculation
$$
 \int_{\tilde{{\mathcal Z}}_{s}}
 \tilde{\alpha}^{n+1}\barpulltwoforms^{n-1}
 = (n+1)\quadrics (\alpha )
   \left(
   \int_{\tilde{{\mathcal Z}_{s}}}
   \tilde{\alpha}\pulltwoform_{s}^{n-1}\barpulltwoforms^{n}
   \right)^{n-1}.
$$
 Let 
 $E$ be an irreducible component of the exceptional locus of 
 $\tilde{{\mathcal Z}}_{s} \to {\mathcal Z}_{s}$.
 Assume $\alpha := \phi_{s}\phi^{-1}_{t}(\twoform_{t})$.
 Then $\tilde{\alpha}^{n+1}|_{{\mathcal U}_{s}}$ = 0 because
 $\alpha |_{{\mathcal U}_{s}}$ is a holomorphic two form on 
 the smooth locus ${\mathcal U}_{t}$ of ${\mathcal Z}_{t}$ and
 ${\mathcal U}_{t}$ is deformation 
 equivalent to ${\mathcal U}_{s}$ by \cite[Theorem 8]{namikawa1}.
 By \cite[Remark (1) (p24)]{namikawa1},
 $\pulltwoform_{s}^{n-1}|_{E}$ = 0. Hence the left hand side equals
 zero. If $t$ is very near to $s$, we obtain 
 $   \int_{\tilde{{\mathcal Z}_{s}}}
   \tilde{\alpha}\pulltwoform_{s}^{n-1}\barpulltwoforms^{n} \ne 0
 $. Thus we obtain 
 $\quadrics (\alpha) = 0$. 
\qed

\noindent
 We continue the proof of Lemma \ref{prep1}. By the above Claim,
 $\phi_{s}$ maps an open set of $\quadricsurface_{0}$
 to an open set of $\quadricsurface_{s}$
 isomorphically. Since both quadric surfaces
 are irreducible, we obtain that 
 $\phi_{s} (\quadricsurface_{0}) = \quadricsurface_{s}$.
 Since $\quadrics (\twoform_{s} + \bar{\twoform}_{s}) = 1$,
 we are done.
\qed

\noindent
(2.3) The following Lemma is the key of the proof of Theorem \ref{main}
 and \ref{main2}
 which is based on arguments in 
 \cite[Lemma 1.9]{bogomolov} and \cite[Theorem 5.6]{Huyhab}.
\begin{lem}\label{key_lemma}
 Let $\nu : \tilde{{\mathcal Z}} \to {\mathcal Z}$ be a resolution
 of ${\mathcal Z}$ and ${\mathcal D}^{\circ}$ the open set
 set of ${\mathcal D}$ such that  the morphism 
 $\tilde{{\mathcal Z}} \to {\mathcal D}$  is smooth over
 ${\mathcal D}^{\circ}$.
 We fix one point $s$ of ${\mathcal D}^{\circ}$.
 Let $\alpha \in H^{2}({\mathcal Z}_{s}, \CC )$  and
 $\tau$  a $(k,k)$-form on $\tilde{{\mathcal Z}}_{s}$. 
 If $\tau$ is a $(k,k)$-form on every 
 fibre $\tilde{{\mathcal Z}}_{t}$ near $s$,
\begin{equation}\label{desired}
   \tau \cdot \nu^{*}\alpha^{2(n-k)} = C_{k} \quadrics (\alpha )^{n-k},
\end{equation}
 where $C_{k}$ is a constant depending on only ${\mathcal Z}_{s}$.
\end{lem}
\noindent
{\scshape Proof. \quad} 
 Let $\cupsurface_s$ be the hypersurface in $H^{2}({\mathcal Z}_{s},\CC )$
 defined by
$$
 \{ \tau \cdot \nu^{*}\alpha^{2(n-k)} = 0 | 
    \alpha \in H^{2}({\mathcal Z}_{s},\CC )\}.
$$
 If we prove that $\quadricsurface_{s} = \cupsurface_{s}$,
 we obtain equation (\ref{desired}) because both hand sides
 of (\ref{desired})
 have same degree and same zero locus.
 We choose an open set $V$ of $\quadricsurface_{s}$
 which is contained in $p_{s}({\mathcal D}^{\circ})$.
 For every point $\alpha \in V$, there exists the point $t$ of 
 ${\mathcal D}^{\circ}$ such that $\alpha$ defines
 a symplectic form on the smooth part of ${\mathcal Z}_{t}$.
 Thus $\tau \cdot (\nu^{*}\alpha )^{n-k+1} = 0$ in 
 $H^{2n+2}(\tilde{{\mathcal Z}}_{t},\CC )$ and hence
 $\tau\cdot (\nu^{*}\alpha )^{n-k+1} = 0$
 in $H^{2n+2}(\tilde{{\mathcal Z}}_{s},\CC )$.
 By analytic continuation, we obtain 
 $\tau\cdot (\nu^{*}\alpha )^{n-k+1}=0$ for every point of 
 $\quadricsurface_{s}$.
 Therefore we obtain $\quadricsurface_{s} \subset \cupsurface_{s}$.
 We prove the opposite inclusion.
 Assume the contrary. Then there exists 
 $\beta \in H^{2}({\mathcal Z}_{s},\CC )$ such that
 $\tau \cdot \nu^{*}\beta^{2(n-k)} = 0$ and 
 $\beta \not\in \quadricsurface_{s}$.
 We choose a general element $\gamma$
 of $H^{2}({\mathcal Z}_{s},\CC )$ 
 such that $\tau \cdot \gamma^{2(n-k)} \ne 0$. Then the line which
 pass $\beta$ and $\gamma$ intersect $\quadricsurface_{s}$
 with two points. Let $\delta_{0}$ and  $\delta_{1}$ be
 these points. We write
$$
 \beta = \lambda_{11}\delta_{0} + \lambda_{12}\delta_{1}, \quad
 \gamma= \lambda_{21}\delta_{0} + \lambda_{22}\delta_{1}.
$$
 We remark that $\lambda_{**} \ne 0$ because
 $\beta$, $\gamma$ and $\delta_{*}$ are mutually distinct.
 Since 
 $\tau\cdot (\nu^{*}\delta_{0})^{n-k+1} 
  = \tau\cdot (\nu^{*}\delta_{1})^{n-k+1} = 0$,
\begin{eqnarray*}
 \tau\cdot (\nu^{*}\beta)^{2(n-k)} 
           &=& {}_{2(n-k)}C_{n-k}\lambda_{11}^{n-k}\lambda_{12}^{n-k} 
                \delta_{0}^{n-k}\delta_{1}^{n-k} = 0 \\
 \tau\cdot (\nu^{*}\gamma)^{2(n-k)}
            &=& {}_{2(n-k)}C_{n-k}\lambda_{21}^{n-k}\lambda_{22}^{n-k}
                \delta_{0}^{n-k}\delta_{1}^{n-k} \ne 0 .
\end{eqnarray*}
 That derives a contradiction.
\qed

\noindent
(2.4){\scshape Proof of Theorem \ref{main}. \quad}
 From Lemma \ref{key_lemma}, we obtain
$$
 C_{Z} \quadrics (\phi_{s}(\alpha))^{n} = \phi_s (\alpha)^{2n}.
$$
 Since $\phi_{s}(\alpha )^{2n} = \alpha^{2n}$,
 we obtain that $C_{Z} \quadric (\alpha )^{n} = \alpha^{2n}$ by
 Lemma \ref{prep1}.
 We investigate the index of $\quadric$. 
 From \cite[Proposition 9]{namikawa1},
 there exists Hodge decomposition of $H^{2}(Z, \CC )$.  
 Let $A$ be an ample divisor on $Z$ and $H$ an element
 of $H^{1,1}(Z, \CC )_{{\mathbb R}}$
 such that $\quadric (H,A) = 0$. We consider the
 following equation:
\begin{equation}\label{nondegenerate}
 C\quadric (\lambda H + A)^{n} = (\lambda H + A)^{2n}.
\end{equation}
 If we compare $\lambda$ term of both hand sides,
 we obtain $H.A^{2n-1} = 0$
 from the assumption $\quadric (A,H) = 0$. By Hodge-Riemann bilinear
 relation,
$$
 H^2 .A^{2n-2} \cdot A^{2n} \le (H.A^{2n-1})^2 = 0.
$$
 Hence $H^2 .A^{2n-2} < 0$ if $H \not\equiv 0$.
 Comparing $\lambda^2$ term of
 the both hand side of the equation (\ref{nondegenerate}),
 we obtain $\quadric (H) < 0$ if $H \not\equiv 0$.
 Combining Lemma \ref{rank}, the index of $\quadric$ is 
 $(3, h^2 (Z, \CC )-3)$.
\qed

\noindent
(2.5) {\scshape Proof of Theorem \ref{main2}. \quad}
 From the proof of (1) of \cite[Theorem 8]{namikawa1},
$$
 \phi_{s} : H^{2}(Z,\CC ) \cong H^{2}({\mathcal Z},\CC ) \to 
            H^{2}({\mathcal Z}_{s},\CC )
$$
 is isomorphism. Hence
 $\chi_{Z} (D) = \chi_{{\mathcal Z}_{s}}(\phi_{s} (D))$.
 By Lemma \ref{prep1},
 it is enough to prove that 
 there exists constants $C_k$ $(0 \le k \le n)$
 and they satisfy
$$
 \chi_{{\mathcal Z}_{s}} (D) = \sum_{i=0}^{n}C_k \quadrics (D)^{k},
$$
 for a Cartier divisor $D$ of ${\mathcal Z}_{s}$.
 Let us consider $\chi_{\tilde{{\mathcal Z}}_{s}}(\nu^{*}D)$.
 By \cite[Theorem 6]{namikawa1},  ${\mathcal Z}_{s}$ has only rational
 singularities. Thus $\chi_{{\mathcal Z}_{s}}(D) = 
 \chi_{\tilde{{\mathcal Z}}_{s}}(\nu^{*}D)$.
 By Serre duality
$$
 \chi_{{\mathcal Z}_{s}} (D) =  \chi_{{\mathcal Z}_{s}} (-D).
$$
 Hence each term of Riemann-Roch formula of
 $\chi_{\tilde{{\mathcal Z}}_{s}}(\nu^{*}D)$
 is expressed 
$$
 (\mbox{Polynomial of Chern classes of $\tilde{{\mathcal Z}}_{s}$})
 \cdot (\nu^{*}D)^{2k}.
$$
 Since polynomials of Chern classes remain $(k,k)$-form
 under small deformation, we obtain
 each term of Riemann-Roch formula is expressed as
 $C_k \quadrics (D)^{k}$ by Lemma \ref{key_lemma}. 
\qed

\section{Fibre space structure}
\noindent
(3.1) We prove Corollary \ref{application} in three steps:
\begin{description}
 \item[Step 1] $\dim B = n$.
 \item[Step 2] For a general fibre $F$ of $f$, the restriction
                of a symplectic form on the smooth locus of $F$ is
                identically zero.
 \item[Step 3] $F$ is an Abelian variety and $f(Z_{sing})$ forms
               a proper closed subset of $B$.
 \item[Step 4] $\rho (B) = 1$
 \item[Step 5] For every effective divisor $D$ on $Z$, $\dim f(D) \ge n-1$.
\end{description}
 We fix some notations. Let $A$ be an ample divisor on $Z$ and $H$
 an ample divisor on $B$.

\noindent
(3.2) {\scshape Step 1. \quad}
 From Theorem \ref{main}, 
$$
 C_{Z}\quadric (A + \lambda f^{*}H)^{n} = (A + \lambda f^{*}H)^{2n}.
$$
 Since $C_{Z}\quadric (f^{*}H) = (f^{*}H)^{2n} = 0$,
 we obtain
\begin{equation}
 C_{Z}(\quadric (A) + 2\lambda \quadric (A,f^{*}H))^{n}
  = (A + \lambda f^{*}H)^{2n}.
\end{equation}
 If we compare both hand sides the above equation,
 we obtain
\begin{eqnarray*}
  A^{k}.(f^{*}H)^{2n-k} &=& 0, \quad (k < n) \\
  2n A^{2n-1}.f^{*}H &=& C_{Z} \quadric (A, f^{*}H ) \cdot\quadric (A)^{n-1}.
\end{eqnarray*}
 By Lemma \ref{rank}, $\quadric (A) > 0$. Hence
 $\quadric (A,f^{*}H) > 0$ and
$$
 A^{k}.(f^{*}H)^{2n-k} > 0, \quad (k \ge n).
$$
 Hence $\dim B = n$.

\noindent
(3.3) {\scshape Step 2. \quad}
 Let $F$ be a general fibre of $f$ and
 $U_{F}$ the smooth part of $F$.
 In order to prove $\twoform |_{U_{F}} \equiv 0$,
 we prove
$$
 \int_{F} \twoform \wedge \bar{\twoform} A^{n-2} = 0.
$$
 We have
$$
 \int_{F} \omega \wedge \bar{\omega} A^{n-2} = 
 c (\omega  \bar{\omega} A^{n-2} (f^{*}H)^n ),
$$
 where $c$ is a nonzero constant.
 Hence we will show $\omega \bar{\omega} A^{n-2} (f^{*}H)^n = 0$.
 By Theorem \ref{main},
$$
 C_{Z} \quadric (\twoform + \bar{\twoform } + sA + tf^{*}H,
   \twoform + \bar{\twoform } + sA + tf^{*}H )^n 
 = ( \twoform +\bar{\twoform } + sA + tf^{*}H )^{2n}.
$$
 Calculating the left hand side, we obtain
$$
 C_{Z} (\quadric (\twoform + \bar{\twoform }) + s^2 \quadric (A)
      + 2s\quadric (\twoform + \bar{\twoform },A) + 
        2t\quadric(\twoform  + \bar{\twoform },f^{*}H)  +
        2st\quadric (A,f^{*}H))^n . 
$$
 From the definition of $\quadric$ in Theorem \ref{main},
$$
 \quadric(\omega + \bar{\omega}, A) = 
 \quadric(\omega + \bar{\omega}, f^{*}H) = 0 .
$$
 Thus
 we conclude that  $\omega \bar{\omega} A^{n-2} (f^{*}H)^n = 0$ by
 comparing the $s^{n-2}t^n$ term of both hands sides.

\noindent
(3.4) {\scshape Step 3. \quad}
 Let $U$ be the smooth locus of $Z$.
 If we choose a point $x$ of $B$ generally,
 $f^{\circ} : U \to B$ is smooth at $x$ and
 $U_{F} = F \cap U$.
 We consider the following diagram:
$$
 \begin{array}{cccccc}
  0 \to & f^{*}T_{B,x} & \to  & T_{U}& \to &  T_{F_{U_{F}}} \to 0 \\
        &              &      & \downarrow & &                    \\
  0 \to &\Omega^{1}_{U_{F}} & \to & 
         \Omega^{1}_{U} & \to & f^{*}\Omega^{1}_{B,x} \to 0.
 \end{array}
$$
 From the above diagram and Step 2, $h^{0}(T_{U_{F}}) = n$. 
 Then $F$ is an Abelian variety by the following Lemma.
\begin{lem}
 Let $F$ be a normal variety such that
 $K_F \sim {\mathcal O}_{F}$ and
 $h^{0}(T_{U_{F}}) = \dim F$, where $U_{F}$ is
 the smooth locus of $F$. Then $F$ is an Abelian variety. 
\end{lem}  
\noindent
{\scshape Proof. \quad}
 Since $F$ is normal,
 $\dim (F\setminus U_{F}) \le \dim F - 2$.
 Hence $h^{0}(\Theta_{F}) = \dim F$ by analytic continuation.
 Let $\tilde{F}$ be a resolution of $F$. Then
 $h^{0}(T_{\tilde{F}}) = h^{0}(\Theta_{F}) = \dim F$.
 Since $K_{F} \sim {\mathcal O}_{F}$,
 we have an injection $T_{\tilde{F}} \to \Omega^{\dim F
 -1}_{\tilde{F}}$.
 Hence $h^{\dim F -1}({\mathcal O}_{\tilde{F}}) = 
 h^{0}(\Omega^{\dim F - 1}_{\tilde{F}}) \ge \dim F$.
 By Serre duality 
 $h^{\dim F -1}({\mathcal O}_{\tilde{F}}) =
  h^{1}(K_{\tilde{F}}) $
 From Grauert-Riemenschneider vanishing theorem,
 $ h^{1}(K_{\tilde{F}}) = 
 h^{1}(\pi_{*}K_{\tilde{F}}) =
 h^{1}({\mathcal O}_{F})$. 
 Hence $h^{1}(F,{\mathcal O}_{F}) \ge \dim F$.
 By \cite[Theorem 13]{Abelian_Characterization},
 $F$ is an abelian variety. 
\qed

\noindent
 Since $F$ is a complete intersection in $Z$,
 $Z$ is smooth in a neighbourhood of $F$ and
 $f(Z_{sing})$ forms a proper closed subset of $B$.

\noindent
(3.5){\scshape Step 4. \quad}
 Let $D$ be a Cartier divisor on $B$.
 We prove that
 $f^{*}(D + \lambda H) \equiv 0$ for a suitable number $\lambda$.
 From Theorem \ref{main}, $\quadric $ is nondegenerate. Hence
 it is enough to prove
 that
$$
 \quadric (f^{*}(D + \lambda H)) = 
 \quadric (f^{*}(D + \lambda H), A) = 0.
$$
 From Theorem \ref{main}, we obtain 
 $C_{Z} \quadric (f^{*}(D + \lambda H))^n =
 (f^{*}(D + \lambda H))^{2n} = 0$ 
 for every $\lambda$. Hence if we choose $\lambda$ suitably,
 we obtain $ \quadric (f^{*}(D + \lambda H), A) = 0$.
\qed

\noindent
(3.6) {\scshape Step 5. \quad}Let $D$ be an effective divisor on $Z$.
 We derive a contradiction assuming that $\dim f(D) < \dim B -1$.
 Under this assumption, we obtain that $D.f^{*}H^{n-1}$ is
 numerically trivial.
 The following equation
$$
 C_{Z}\quadric (D + tf^{*}H)^n = (D + tf^{*}H)
$$
 tells us that $\quadric (D,f^{*}H) = 0$ because
 $D^{n}f^{*}H^{n} = 0$.
 We consider the following equation:
$$
 (s D + t f^{*}H + A)^{2n} = 
 C_{Z}\quadric (s D + t f^{*}H + A)^{n}.
$$
 Comparing $st^{n-1}$-term of the both hand sides, we obtain that 
 $\quadric (D,A) \quadric (A,f^{*}H) = 0$.
 For an effective divisor $E$ on $Z$, we consider the following
 equation:
$$
 C_{Z}\quadric (tE + A)^{n} = 
 (tE + A)^{2n}.
$$
 Comparing $t$-term of the both hand sides,
 we obtain that $\quadric (D,A) > 0$ and $\quadric (f^{*}H,A) > 0$.
 That derives a contradiction.
\qed

\noindent
(3.7) Under the conditions of Corollary \ref{application},
  there exists an example such that $-K_B$ is not ample.

\noindent
{\scshape Example. \quad}
 Let $E_{i}$ be an elliptic curve. We consider the abelian
 6-fold 
 $\tilde{Z} := E_{1}\times \cdots \times E_{6}$. 
 Let $G$ be the finite subgroup
 of ${\mathrm G}{\mathrm L}(6,\CC )$ generated by
$$
\left(
\begin{array}{cccccc}
 \zeta_{6} & 0 & 0 & 0 & 0 & 0 \\
 0  &\zeta^{5}_{6} & 0 & 0 & 0 & 0 \\
 0 & 0 & -1 & 0 & 0 & 0 \\
 0 & 0 & 0 & -1 & 0 & 0 \\
 0 & 0 & 0 & 0 & \zeta_{3} & 0 \\      
 0 & 0 & 0 & 0 & 0 & \zeta^{2}_{3} 
\end{array} 
\right) \mbox{and}
\left(
\begin{array}{cccccc}
 0 & 0 & 1 & 0 & 0 & 0 \\
 0 & 0 & 0 & 1 & 0 & 0 \\
 0 & 0 & 0 & 0 & 1 & 0 \\
 0 & 0 & 0 & 0 & 0 & 1 \\
 1 & 0 & 0 & 0 & 0 & 0 \\      
 0 & 1 & 0 & 0 & 0 & 0 
\end{array} 
\right),
$$
 where $\zeta_{n}$ means a $n$-th root of unity.
 Then the quotient $Z := \tilde{Z}/G$ satisfies
 the conditions Theorem \ref{main} and it admits
 a fibration $Z \to (E_1 \times E_3 \times E_5 )/G$. By direct calculation,
 $K_{E^{3}/G} \equiv 0$.

\end{document}